\documentclass[12pt,a4paper]{amsart}
\usepackage{amsmath}

\newcommand{\fsl}{\mathfrak{sl}}
\newcommand{\fso}{\mathfrak{so}}

\newcommand{\bZ}{\mathbb{Z}}
\newcommand{\bQ}{\mathbb{Q}}

\newcommand{\cA}{\mathcal{A}}
\newcommand{\cB}{\mathcal{B}}
\newcommand{\cT}{\mathcal{T}}

\newtheorem{theorem}{Theorem}[section]

\newtheorem{defn}[theorem]{Definition}

\DeclareMathOperator{\Hom}{Hom}
\DeclareMathOperator{\End}{End}

\begin{document}
\title{Invariant tensors and cellular categories}
\author{Bruce W. Westbury}
\address{Mathematics Institute\\
University of Warwick\\
Coventry CV4 7AL}
\email{Bruce.Westbury@warwick.ac.uk}
\date{\today}
\begin{abstract} Let $U$ be the quantised enveloping algebra associated
to a Cartan matrix of finite type. Let $W$ be the tensor product of
a finite list of highest weight representations of $U$. Then $\End_U(W)$
has a basis called the dual canonical basis and this gives an integral form for
 $\End_U(W)$. We show that this integral form is cellular by using results
due to Lusztig.
\end{abstract}
\maketitle

\section{Introduction}
The definition of a cellular algebra was given in \cite{MR1376244}
and the examples of cellular algebras given were the Temperley-Lieb
algebras, the Hecke algebras, the Brauer (or Birman-Wenzl) algebras
and the Ariki-Koike algebras.

The aim of this paper is to give a general construction of cellular
algebras using the dual canonical bases of spaces of invariant tensors
defined in \cite[Part IV]{MR1227098}. This construction includes the Temperley-Lieb
algebras, the Hecke algebras, and the  Birman-Wenzl algebras.
This should not be confused with the observation made in \cite{MR1713038}
that the results of \cite[Chapter 29]{MR1227098} imply that the
algebra $\mathbf{\dot U}$ and its finite dimensional quotients
$\mathbf{\dot U}/\mathbf{\dot U}[P]$ are cellular.

The definition of a cellular algebra has a natural generalisation
to cellular categories. An endomorphism algebra in a cellular category
is cellular and many of the sequences of cellular algebras arise this
way. Examples are the many diagram algebras such as the partition algebras,
the blob algebras, the affine Temperley-Lieb algebras, the contour algebras,
the two-boundary Temperley-Lieb algebras.

Let $U$ be the quantised enveloping algebra associated
to a Cartan matrix of finite type and let $V$ be a self-dual finite
dimensional highest weight representation. Then the category of invariant
tensors has objects $n=0,1,2,\ldots $ and the morphisms are
$\Hom_U(\otimes^nV,\otimes^mV)$. The examples which have been studied
are
\begin{itemize}
\item \cite{MR1446615}, the representations of $\fsl(2)$ 
\item \cite{MR992598}, the vector representations of the symplectic and orthogonal Lie algebras
\item \cite{MR1403861}, the fundamental representations of the rank two simple Lie algebras
\item \cite{MR1941983}, the adjoint representation of $\fsl(n)$
\item \cite{MR2388243}, the spin representation of $\fso(7)$
\end{itemize}
In each case there is an integral form of the category of invariant
tensors which is constructed using planar graphs and in each case
this integral form can be seen to be cellular.

If $V$ is not self-dual then we need to consider tensor products
of copies of both $V$ and $V^*$. The only case that has been studied
is for $V$ the vector representation of $\fsl(n)$. Then the walled
Brauer algebras are endomorphism algebras in the category of invariant
tensors.



In this paper we consider a larger category of invariant tensors.
This is the full subcategory of the category of finite dimensional
representations of $U$ whose objects are tensors products of highest
weight representations. In this paper we observe that the results of
\cite[Part IV]{MR1227098} show that the dual canonical basis defines an
integral form of this category and that this integral form is cellular.
In particular this shows that if $W$ is any tensor product of highest
weight representations then the integral form of $\End_U(W)$ associated
to the canonical basis is cellular.

This work was influenced by the principle that the diagram bases and the
dual canonical bases have similar properties. In the examples in which
the diagram bases have been defined it is straightforward to see that
the category is cellular. It is then natural to look for a cellular
structure in the dual canonical basis.

The diagram bases and dual canonical bases for $\fsl(2)$ are discussed and
shown to be the same in \cite{MR1446615}. The diagram
bases and the dual canonical bases for the category of
invariant tensors of $\fsl(3)$ are compared in \cite{MR1680395} and shown to be
different.

The diagram bases and the dual canonical bases have the following properties in
common. They each define an integral form of the category of invariant tensors.
They each admit a cellular structure with the same poset. In each case the
tensor functor $-\otimes -$ and the duality functor respect the cell structure.
This makes it reasonable to conjecture that these categories are equivalent
as integral forms and that the equivalences respect the cellular structure, the
tensor structure and the duality.

\section{Cellular categories}
First we extend the definition of a cellular algebra given in
\cite{MR1376244} to cellular categories. The motivating examples
are the Temperley-Lieb category discussed in \cite{MR1343661},
and the categories whose endomorphism algebras are the sequences
of algebras in the Introduction.

\begin{defn}
Let $R$ be a commutative ring with identity. Let $\cA$ be a $R$-linear
category with an anti-involution $*$. Then cell datum for $\cA$ consists
of a partially ordered set $\Lambda$, a finite set $K(n,\lambda)$ 
for each $\lambda\in\Lambda$ and each object $n$ of $\cA$, and for 
$\lambda\in\Lambda$ and $n$,$m$ any two objects of $\cA$ we have an inclusion
\[ C\colon K(n,\lambda)\times K(m,\lambda)\rightarrow \Hom_\cA(n,m) \]
\[ C\colon (S,T)\mapsto C^\lambda_{S,T} \]
The conditions that this datum is required to satisfy are:
\begin{enumerate}
\item[C-1] For all objects $n$ and $m$, the image of the map
\[ C\colon \coprod_{\lambda\in\Lambda} K(n,\lambda)\times K(m,\lambda) 
\rightarrow \Hom_\cA(n,m) \]
is a basis for $\Hom_\cA(n,m)$ as an $R$-module.
\item[C-2] For all objects $n,m$, all $\lambda\in\Lambda$ and all
$S\in K(n,\lambda)$, $T\in K(m,\lambda)$ we have
\[ \left(C^\lambda_{S,T}\right)^*=C^{\lambda}_{T,S} \]
\item[C-3] For all objects $p,n,m$, all $\lambda\in\Lambda$ and all
$a\in \Hom_\cA(p,n)$,
$S\in K(n,\lambda)$, $T\in K(m,\lambda)$ we have
\[ aC^\lambda(S,T)=\sum_{S^\prime\in K(p,\lambda)}
r_a(S^\prime,S)C^\lambda_{S^\prime ,T}
\mod \cA(<\lambda) \]
where $r_a(S^\prime,S)\in R$ is independent of $T$ and 
$\cA(<\lambda)$ is the $R$-linear span of
\[ \left\{ C^\mu_{S,T} | \mu < \lambda ; S\in K(p,\mu),T\in K(m,\mu)
\right\} \]
\end{enumerate}
\end{defn}
Another way of formulating C-3 is that for all objects
$p,n,m$, all $\lambda\in\Lambda$ and all
$S\in K(p,\lambda)$,$T\in K(n,\lambda)$,
$U\in K(n,\lambda)$,$V\in K(m,\lambda)$ we have
\[ C^\lambda(S,T)C^\lambda(U,V)=\langle T,U\rangle C^\lambda_{S,V}
\mod \cA(<\lambda) \]
where $\langle T,U\rangle$ is independent of $S$ and $V$.

A consequence of this definition is that for all order ideals
$\pi\subset\Lambda$ we have an ideal $\cA(\pi)$ where
for all objects $n,m$ the subspace $\cA(\pi)(n,m)$ is the
$R$-linear span of 
\[ \left\{ C^\mu_{S,T} | \mu\in\pi ; S\in K(n,\mu),T\in K(m,\mu) \right\} \]

This is a generalisation of the definition of a cellular algebra
since we can regard any algebra over $R$ as an $R$-linear category
with one object.

More significantly, if $\cA$ is a cellular category then $\End(n)$
is a cellular algebra for any object $n$ of $\cA$.

A functor $\phi\colon\cA\rightarrow\cA^\prime$ between cellular categories
is cellular if we have a map of partially ordered sets
$\phi\colon \Lambda\rightarrow\Lambda^\prime$ such that for any order
ideal $\pi\in\Lambda$ we have $\phi(O(\pi))\subset O(\phi(\pi))$.

One reason for considering cellular categories instead of just 
the cellular endomorphism algebras is that for each object $n$
and each order ideal $\pi\subset\Lambda$ we have an $R$-linear functor
$\rho(n;\pi)$ from $\cA$ to the category of left
$\End_\cA(n)$-modules which on objects is given by
\[ \rho(n;\pi)\colon m\mapsto \Hom_{\cA(\pi)}(n,m) \]
Furthermore if $\pi\subset\pi^\prime$ then we have a natural transformation
from $\rho(n;\pi)$ to $\rho(n;\pi^\prime)$. This is used
in \cite{MR1659204}.

\section{Canonical bases}
The rest of the paper is based on \cite[Part IV]{MR1227098}.
In particular all references are to this book and we have tried
to use notation which is the same as or compatible with the
notation in this book.

Let $C$ be a Cartan matrix of finite type and let $U$ be the 
Drinfeld-Jimbo quantum group associated to $C$. This is a Hopf
algebra over the field $\bQ(v)$. The integral weight lattice is denoted
by $X$ and the dominant integral weights are denoted by $X^+$.
For each $\lambda\in X^+$ there is a finite dimensional representation
of $U$ with highest weight $\lambda$; we denote this representation
by $V_\lambda$. This representation has a canonical basis which we denote
by $B_\lambda$. There is an anti-involution, $\omega$, of the
$\bQ(v)$-algebra $U$. This induces a duality on representations
which we denote by $M\mapsto {}^\omega M$. Then
${}^\omega V_\lambda\cong V_{-w_o\lambda}$ where $w_0$ is the
longest element of the Weyl group so we put ${}^\omega\lambda = -w_0\lambda$.

A based module is defined in \cite[27.1.2]{MR1227098}. Here we
consider a category $\cB$. The objects are the based modules.
Let $(M,B)$ and $(M^\prime,B^\prime)$ be based modules. Then a
morphism $(M,B)\rightarrow (M^\prime,B^\prime)$ is a homomorphism
of $U$-modules $\phi\colon M\rightarrow M^\prime$ which maps
$B$ to the $\bZ[v,v^{-1}]$-module generated by $B^\prime$ and
which is compatible with the involutions in the sense that
$\overline{\phi(m)}=\phi(\overline{m})$ for all $m\in M$.

Put $\delta=v+v^{-1}$. Then observe that $\bZ[\delta]$ is the
subring of $\bZ[v,v^{-1}]$ fixed by the bar involution. Then
it is clear that $\cB$ is a $\bZ[\delta]$-linear category.
Then \cite[27.3]{MR1227098} constructs a tensor product of based modules
and shows that this makes $\cB$ a monoidal category.
Let $(M,B)$ and $(M^\prime,B^\prime)$ be based modules. Then the
tensor product is the
based module $(M\otimes M^\prime, B\diamond B^\prime)$.

For each $\lambda\in X^+$ we have a based module $(V_\lambda ,B_\lambda)$
where $B_\lambda$ is the canonical basis of $V_\lambda$. Furthermore
\cite[Proposition 21.1.2]{MR1227098} says that there is also a based module
$({}^\omega V_\lambda,{}^\omega B_\lambda)$. Let $\cT$ be the full 
subcategory of $\cB$ whose objects are finite tensor products of based
modules isomorphic to these based modules. Then, by construction, $\cT$ is  a
$\bZ[\delta]$-linear monoidal subcategory of $\cB$. This category is also closed under
the duality $\omega$ and is in fact a spherical category, as defined in
\cite{MR1686423}.

Then we show that the $\bZ[\delta]$-linear category $\cT$ with the
anti-involution $\omega$ is a cellular category. The poset $(\Lambda,<)$
is the poset $(X^+,>)$.

If $(M,B)$ is a based module \cite[27.2.1]{MR1227098}
gives a partition of $B$ into subsets $B[\lambda]$.
Furthermore \cite[27.2.5]{MR1227098} shows that $B[0]$ can be regarded
as a base for the space of coinvariants of $M$.
There is a canonical identification of $\Hom_U(M,M^\prime)$ with the
dual of the space of coinvariants of $M\otimes {}^\omega M^\prime$.
Hence we may regard $(B\diamond {}^\omega B^\prime)[0]$ as a basis
of $\Hom(M,M^\prime)$.

The subsets $B[\lambda]^{hi}\subset B[\lambda]$ and $B[\lambda]^{lo}\subset B[\lambda]$
are defined in \cite[27.2.3]{MR1227098}
and \cite[Proposition 27.2.6]{MR1227098} says that for $\lambda=0$
both of these inclusions are bijections.

Let $(M,B)$  be an object of
$\cT$. Then we define the cells by
\[ K((M,B),\lambda)=B[\lambda]^{hi} \]

The map $C$ is defined by
\[ C\colon (b,b^\prime)\mapsto  b\diamond {}^\omega b^\prime \]

This gives cell datum and it remains to check that this satisfies the conditions for
a cellular category. The condition C-2 holds by construction.

The partition of $(B\diamond B^\prime)[0]$ as 
\[ (B\diamond B^\prime)[0] = \cup_{\mu\in X^+}
\{ b\diamond b^\prime |
b\in B[-w_0(\mu)]^{lo},b^\prime\in B^\prime[\mu]^{hi} \} \]
is given in \cite[Proposition 27.3.8]{MR1227098}.
This shows that condition C-1 holds.

It remains to discuss C-3. This condition holds modulo an ideal $I$.
This means that, for any $M$, that we can replace $M$ by $M/MI$ and
so we assume that $MI=0$. Equivalently we assume that $M[>\lambda]=0$.

By \cite[27.2.4]{MR1227098}, each $b\in B[\lambda]^{hi}$
gives a homomorphism of based modules $\left| b\right\rangle
\colon (M,B)\rightarrow (V_\lambda,B_\lambda)$. Then
$b\diamond {}^\omega b^\prime\colon (M,B)\rightarrow (M^\prime,B^\prime)$
is the composite of the two maps
\[ \left| b\right\rangle
\colon (M,B)\rightarrow (V_\lambda,B_\lambda) \qquad
\left\langle {}^\omega b^\prime\right| \colon (V_\lambda,B_\lambda)
\rightarrow (M^\prime ,B^\prime ) \]
and so we have $b\diamond {}^\omega b^\prime=
\left| b\right\rangle\left\langle {}^\omega b^\prime\right|$.

Now $\End((V_\lambda,B_\lambda))$ is a free $\bZ[\delta]$-module
with basis the identity map $1_\lambda$. 


\end{document}